\documentclass[11pt, reqno]{amsart}
\usepackage{enumitem}
\usepackage{amsthm,amsmath,amssymb}

\usepackage{mathrsfs}

\numberwithin{equation}{section}

\newtheorem{theorem}{Theorem}[section]

\theoremstyle{plain}
\newtheorem{proposition}[theorem]{Proposition}

\newtheorem{corollary}[theorem]{Corollary}
\newtheorem{remark}[theorem]{Remark}

\def\be{\begin{equation}}
	\def\ee{\end{equation}}

\def\d{\nabla}

\DeclareMathOperator{\Ric}{Ric}

\DeclareMathOperator{\vol}{Vol}

\DeclareMathOperator{\Hess}{Hess}

\numberwithin{equation}{section}

\begin{document}
\title[]{Time analyticity for the heat equation under  Bakry-\'Emery Ricci curvature condition}
\author{Ling Wu}
\address{School of Mathematical Sciences and Shanghai Key Laboratory of PMMP, East China Normal University, Shanghai 200241, China}
\email{ 52215500012@stu.ecnu.edu.cn}
\date{}

\begin{abstract}
	Inspired by Hongjie Dong and Qi S. Zhang's article \cite{ZQ2}, we find that the analyticity in time for a smooth solution of the heat equation with exponential quadratic growth in the space variable can be extended to any complete noncompact Riemannian manifolds with Bakry-\'Emery Ricci curvature bounded below and the potential function being of at most quadratic growth. Therefore, our result holds on all gradient Ricci solitons. As a corollary, we give a necessary and sufficient condition on the solvability of the backward heat equation in a class of functions with the similar growth condition. In addition, we also consider the solution in certain $L^p$ spaces with $p\in[2,+\infty)$ and prove its analyticity with respect to time.

\end{abstract}
\maketitle

\section{Introduction}
Let $(M^{n},g)$ be an $n$-dimensional Riemmanian manifold. The Bakry-\'Emery Ricci curvature tensor of $M$ (\cite{BE}) is defined as
\be\label{BE curvature}
\Ric_f:=\Ric+\Hess f,
\ee
where $f$ is a smooth function on $M$ (called the potential function), and $\Ric$ and $\Hess f$ denote the Ricci curvature tensor  and the hessian of $f$, respectively. It is clear that when $f$ is a constant, $\Ric_f$ reduces to the Ricci curvature tensor. A gradient Ricci soliton is a Riemannian manifold $(M^n, g)$ with constant Bakry-\'Emery Ricci curvature, namely,
\be\label{gradient Ricci soliton}
\Ric + \Hess f=\lambda g
\ee
for some constant $\lambda$. It is called a shrinking, steady, or expanding Ricci soliton when $\lambda>0$, $=0$, or $<0$, respectively. Also, manifolds with lower Bakry-\'Emery Ricci curvature bound are closely related to the singularity analysis of the Ricci flow and Ricci limit spaces (see e.g., \cite{Ham, Per, St1, St2, LV}). Therefore, many efforts have been made in extending the results under Ricci curvature condition to Bakry-\'Emery Ricci curvature condition.

  The study of the analyticity of the heat equation has a rich history. For generic solutions, as expected, the space analyticity is valid. However, the time analyticity is more delicate and is indeed invalid, since in the Euclidean space, it is easy to construct a non-time-analytic solution of the heat equation in a finite space-time cylinder. Therefore, it is meaningful to study the time analyticity of the heat equation.

Recently, Qi S. Zhang \cite{ZQ1} discovered for any complete noncompact Riemannian whose Ricci curvature is bounded from below, any solution to the heat equation with exponential growth in the space variable is analytic in time. This result was improved to any solution with exponential quadratic growth by Hongjie Dong and Qi S. Zhang \cite{ZQ2}. In particular, as a corollary to this result, they gave a sufficient and necessary condition on the solvability of the backward heat equation. In \cite{WJY}, Jiayong Wu obtained a similar result on the time analyticity of the heat equation on the complete noncompact gradient shrinking Ricci solitons. For more results, see \cite{Gi}, \cite{LSC}, \cite{WLL}, \cite{WD} and references therein.

In \cite{ZQ2}, for Riemannian manifolds with Ricci curvature bounded below the key estimate for the proof of time analyticity of the heat equation is the parabolic mean value inequality, which can also be found in \cite{SWZ} under Bakry-\'Emery Ricci curvature condition. Here we emphasize our result generalizes Hongjie Dong and Qi S. Zhang's result \cite{ZQ2} and can be extend to all gradient Ricci solitons.

\begin{theorem}\label{1.1}
		Let $(M^n,g)$ be a complete noncompact Riemannian manifold with $\Ric_f\ge-Kg$ for some constant $K\ge0$. For a fixed point $o\in M$, assume that there exist non-negative constants $a$ and $b$ such that 
	\begin{equation}\label{sss1}
		|f(x)|\le ad^2(x,o)+b\ for\ all\ x\in M,
	\end{equation}
where $d(x,o)$ is the distance function from $x$ to $o$. Let $u(x,t)$ be a smooth solution of the heat equation $(\Delta-\partial_t)u=0$ on $M\times[-2,0]$ and satisfy exponential quadratic  growth in  the space variable, i.e.,
\begin{equation}\label{s2}
	|u(x,t)|\le A_1e^{A_2d^2(x,o)}\ for\ all\ (x,t)\in M\times[-2,0],
\end{equation}
 where $A_1$ and $A_2$ are some positive constants. Then $u(x,t)$ is analytic in time $t\in[-1,0]$ with radius $\delta>0$ depending only on $n,K,a,b$ and $A_2$. Besides, we have  
\begin{equation}\label{3}
	u(x,t)=\sum_{j=0}^{\infty}a_j(x)\frac{t^j}{j!}
\end{equation} 
with $\Delta a_j(x)=a_{j+1}(x)$ and 
\begin{equation}\label{www1}
	|a_j(x)|\le A_1A_3^{j+1} (j+1)^je^{A_4d^2(x,o)},\ j=0,1,2,...,
\end{equation}
where $A_3$ and $A_4$ are two positive constants depending on $K,n,a,b,A_2$ and $n,a,A_2$, respectively.
\end{theorem}
\begin{remark}
If the potential function is $0$, i.e., $a=b=0$ in \eqref{sss1}, after careful calculation, then we get $A_4=2A_2$ in \eqref{www1}. Theorem \ref{1.1} reduces Honejie Dong and Qi S. Zhang's result \cite{ZQ2}. 
\end{remark}
\begin{remark}
	The growth condition \eqref{s2} is sharp due to the Tychonov's solution of the heat equation in $R^n\times(-\infty,+\infty)$ (see Remark 2.3 in \cite{ZQ2}).
\end{remark}
The conditions in the above theorem are especially satisfied on gradient Ricci solitons. 
For gradient Ricci solitons, it is well known that
\be\label{w1}
S+|\d f|^2=2\lambda f + C,
\ee
where $S$ is the scalar curvature of $M$, $\nabla f$ is the gradient of $f$ and $C$ is a constant. \\
For gradient shrinking solitons, it is showed in \cite{Chen} that $S\ge0$, then setting $\tilde{f}=f+\frac{C}{2\lambda}$, \eqref{w1} implies 
\begin{equation}
	|\nabla \tilde{f}|^2\le 2\lambda \tilde{f},\nonumber
\end{equation}
so   
\begin{equation}
 |\tilde{f}(x)|\le \lambda d^2(x,o)+2|\tilde{f}(o)|.\nonumber
\end{equation}
For gradient expanding solitons, it is showed in \cite{PRS,ZS} that $S\ge n\lambda$, then \eqref{w1} implies that 
\begin{equation}
	\bigg|\nabla \sqrt{-f-\frac{C-n\lambda}{2\lambda}}\bigg|\le \sqrt{-\frac{\lambda}{2}}.\nonumber
\end{equation}
Hence
\begin{equation}
	\sqrt{-f(x)-\frac{C-n\lambda}{2\lambda}} \le \sqrt{-\frac{\lambda}{2}}d(x,o)+\sqrt{-f(o)-\frac{C-n\lambda}{2\lambda}},\nonumber 
\end{equation}
setting $\tilde{f}=f+\frac{C-n\lambda}{2\lambda}$, then 
\begin{equation}
	|\tilde{f}(x)|\le-\lambda d^2(x,o)+2|\tilde{f}(o)|.\nonumber
\end{equation}
For gradient steady solitons, we know $S\ge0$ in \cite{Chen}.\\
If $C=0$ in \eqref{w1}, then $f$ is a constant.\\
If $C\neq 0$ in \eqref{w1}, by scaling the metric $g$, we can get  
\begin{equation}
	S+|\nabla f|^2=1, \nonumber
\end{equation}
which implies
\begin{equation}
	|f(x)|\le d(x,o)+|f(o)|\nonumber.
\end{equation}
To sum up, for gradient Ricci solitons, we can always adjust $f$ or $g$ to make 
\begin{equation}
	|f(x)|\le a d^2(x,o)+b,
\end{equation}
where $a$ and $b$ are two positive constants depending on $\lambda$ and $f(o)$, respectively.

Therefore, Theorem \ref{1.1} implies the analyticity in time for smooth solutions of the heat equation on  complete noncompact gradient Ricci solitons.

\begin{theorem}\label{s1.3}
	Let $(M^n,g)$ be a complete noncompact gradient Ricci soliton satisfying \eqref{gradient Ricci soliton}. Let $u(x,t)$ be a smooth solution of the heat equation $(\Delta-\partial_t)u=0$ on $M\times[-2,0]$ and satisfy the growth condition
\begin{equation}\label{2}
	|u(x,t)|\le A_1e^{A_2d^2(x,o)}\ for\ all\ (x,t)\in M\times[-2,0],
\end{equation}
where $A_1$ and $A_2$ are some positive constants, and $d(x,o)$ is the distance function from $x$ to a fixed point $o$. Then $u(x,t)$ is analytic in time $t\in[-1,0]$ with radius $\delta>0$ depending only on $n,\lambda,f(o)$ and $A_2$. Besides, we have  
\begin{equation}\label{ww3}
	u(x,t)=\sum_{j=0}^{\infty}a_j(x)\frac{t^j}{j!}
\end{equation} 
with $\Delta a_j(x)=a_{j+1}(x)$ and 
\begin{equation}\label{ww2}
	|a_j(x)|\le A_1A_3^{j+1} (j+1)^je^{A_4d^2(x,o)},\ j=0,1,2,...,
\end{equation}
where $A_3$ and $A_4$ are two positive constants depending on $n,\lambda,f(o),A_2$ and $n,\lambda,A_2$, respectively.
\end{theorem}
\begin{remark}
	In \cite{WJY}, Jiayong Wu obtained a similar result on the time analyticity of the heat equation on the complete noncompact gradient shrinking Ricci solitons. More precisely, he showed that the bound of $a_j(x)$ in \eqref{ww3} is
	\begin{equation}\label{ww4}
		|a_j(x)|\le A_1e^{-\frac{\mu}{2}}e^{\frac{f(x)}{2}}(f(x)+1)^{\frac{n}{4}}A_3^{j+1}j^je^{2A_2d^2(x,o)},\ j=0,1,2,...,
	\end{equation}
where $A_3$ is a constant depending on $n$ and $A_2$ and $\mu=\mu(g,1)$ denotes Perelman's entropy functional.

Comparing \eqref{ww2} with \eqref{ww4}, it is not difficult to find that our result does not depend on $\mu$.
\end{remark}
As an application to Theorem \ref{1.1}, we give a solvable result to the backward heat equation.
 \begin{corollary}\label{1.2}
 	Let $(M^n,g)$ be a complete noncompact Riemannian manifold with $\Ric_f\ge-Kg$ for some constant $K\ge0$. For a fixed point $o\in M$, assume that there exist non-negative constants $a$ and $b$ such that 
 	\begin{equation}\label{sss1}
 		|f(x)|\le ad^2(x,o)+b\ for\ all\ x\in M,
 	\end{equation}
 	where $d(x,o)$ is the distance function from $x$ to $o$. The Cauchy problem for the backward heat equation
 \begin{equation}\label{4}
 	\begin{cases}
 		(\Delta+\partial_t)u=0,& \\
 		u(x,0)=a(x)&
 	\end{cases}
 \end{equation}
 	has a smooth solution with exponential quadratic  growth of the space variable in $M\times(0,\delta)$ for some $\delta>0$ if and only if  
 	\begin{equation}\label{5}
 		|\Delta^ja(x)|\le A_3^{j+1} (j+1)^je^{A_4d^2(x,o)},\ j=0,1,2,...,
 	\end{equation}
 	where $A_3$ and $A_4$ are some positive constants. 
 \end{corollary}
In addition, we also consider the solution of the heat equation in $L^p$ spaces with $p\in[2,+\infty)$ and prove its analyticity with respect to time.

\begin{theorem}\label{1.3}
		Let $(M^n,g)$ be a complete noncompact Riemannian manifold with $\Ric_f\ge-Kg$ for some constant $K\ge0$. For a fixed point $o\in M$, assume that there exist non-negative constants $a$ and $b$ such that 
	\begin{equation}\label{1}
		|f(x)|\le ad^2(x,o)+b\ for\ all\ x\in M,
	\end{equation}
	where $d(x,o)$ is the distance function from $x$ to $o$. Let $u(x,t)$ be a smooth solution of the heat equation $(\Delta-\partial_t)u=0$ on $M\times[-2,0]$. For any $p\ge2$, assume that there exists a positive constant $L$ such that
	\begin{equation}\label{ssss4}
		\left(\int_M |u(x,t)|^pdv\right)^{\frac{1}{p}}\le L\ for\ all\  t\in[-2,0].
	\end{equation}
Then $u(x,t)$ is analytic in time $t\in[-1,0]$ with radius $\delta>0$ depending only on $n,K,a,b$ and $p$. \\
Moreover, we have  
	\begin{equation}\label{3}
		u(x,t)=\sum_{j=0}^{\infty}a_j(x)\frac{t^j}{j!}
	\end{equation} 
	with $\Delta a_j(x)=a_{j+1}(x)$ and 
	\begin{equation}
		|a_j(x)|\le A_6^{j+1} (j+1)^je^{A_7d^2(x,o)}\vol(B_o(1))^{-\frac{1}{p}}L,\ j=0,1,2,...,
	\end{equation}
	where $A_6$ and $A_7$ are two positive constants depending on $n,K,a,b,p$ and $n,a,K,p$, respectively.
\end{theorem}
The rest of this paper is organized as follows. In section 2, we recall a volume comparison theorem and a parabolic mean value inequality from \cite{SWZ} for complete Riemannian manifolds with Bakry-\'Emery Ricci curvature bounded below and the potential function bounded locally. In section 3, applying Hongjie Dong and Qi S. Zhang's method of proof \cite{ZQ2}, we utilize the mean value inequality of section 2 to get a proof of Theorem 1.1, Corollary 1.6 and Theorem 1.7.

\section{Preliminaries}
For a fixed point $o\in M$ and $R>0$, we define 
\begin{equation}
	L(R)=\sup\limits_{B_{o}(3R)}|f|, \nonumber
\end{equation}
where $B_o(3R)$ is the geodesic ball centered at $o\in M$ with radius $3R$.
\begin{theorem}\label{volume element comparison}\cite{SWZ}
	Let $(M^n, g)$ be a complete Riemannian manifold with $\Ric_f\ge-Kg$ for some constant $K\ge0$. Then the following conclusions are true.
	\\
	(a)(Laplacian comparison) Let $r=d(y,p)$ be the distance from any point $y$ to some fixed point $p\in B_{o}(R)$ with $0<r<R$. Then for $0<r_1<r_2<R$, we have
	\begin{equation}\label{Laplacian}
		\int_{r_1}^{r_2}(\Delta r-\frac{n-1}{r})dr \le \frac{K}{6}(r_2^2-r_1^2)+6L(R).
	\end{equation} 
	\\
	(b)(Volume element comparison)Take any point $p\in B_{o}(R)$ and denote the volume form in geodesic polar coordinates centered at $p$ with $J(r,\theta,p)drd\theta$, where $r>0$ and $\theta\in S_p(M)$, a unit tangent vector at $p$. Then for $0<r_1<r_2<R$, we have
	\begin{equation}\label{AC}
		\frac{J(r_2,\theta,p)}{J(r_1,\theta,p)}\le \left(\frac{r_2}{r_1}\right)^{n-1} e^{\frac{K}{6}(r_2^2-r_1^2)+6L(R)} .
	\end{equation}
	\\
	(c)(Volume comparison)For any $p\in B_{o}(R),\ 0<r_1<r_2<R$, we have
	\begin{equation} \label{VC}
		\frac{\vol(B_{p}(r_2))}{\vol(B_{p}(r_1))}\le \left(\frac{r_2}{r_1}\right)^{n} e^{\frac{K}{6}(r_2^2-r_1^2)+6L(R)},
	\end{equation}
	where $\vol(.)$ denotes the volume of a region.
\end{theorem}
From \cite{SWZ}, we have the following parabolic mean value inequality, which is crucial to prove the analyticity of time.

\begin{proposition}(Mean value inequality \cite{SWZ})\label{mean value inequality}
	Let $(M^n,g)$ be a complete Riemannian manifold with $\Ric_f\ge-Kg$ for some constant $K\ge0$. For any real number $s$ and any $0<\tilde{\delta}<\delta'\le1$, let $u$ be a smooth non-negative  subsolution of the heat equation in the cylinder $Q=B_{o}(r)\times(s-r^2,s)$, $0<r<R$.
	
	For $2\le p <\infty$, there exist constants $\tilde{c}_1(n)$ and $\tilde{c}_2(n)$ such that
	\begin{equation}\label{2.4}
		\sup\limits_{Q_{\tilde{\delta}}}u^p\le \frac{\tilde{c}_1(n)e^{\tilde{c}_2(n)(Kr^2+L(R))}}{(\delta'-\tilde{\delta})^{4n}r^2\vol (B_{o}(r))} \cdot\int_{Q_{\delta'}} u^pdvdt.
	\end{equation}
	
	For $0<p<2$, there exist constants $\tilde{c}_3(n,p)$ and $\tilde{c}_4(n)$ such that
	\begin{equation}\label{2.5}
		\sup\limits_{Q_{\tilde{\delta}}}u^p\le \frac{\tilde{c}_3(n,p)e^{\tilde{c}_4(n)(Kr^2+L(R))}}{(\delta'-\tilde{\delta})^{4n}r^2\vol (B_{o}(r))} \cdot\int_{Q_{\delta'}} u^pdvdt.
	\end{equation}
	Here $Q_{\tilde{\delta}}=B_{o}(\tilde{\delta} r)\times(s-\tilde{\delta} r^2,s)$, $Q_{\delta'}=B_{o}(\delta' r)\times(s-\delta' r^2,s)$.
\end{proposition}
\section{proof of the main results}
In this section, we  apply the volume comparison theorem and the parabolic mean value inequality in section 2 to prove the results of this article. We first prove Theorem \ref{1.1}.\\

\noindent{\it Proof of Theorem \ref{1.1}.} Since the heat equation is linear, we can assume that $A_1=1$. Indeed, we just need to prove the time analyticity result at $(x,0)$ for any $x\in M$.

Given $R\ge1$. For any point $x\in B_o(R)$ and a positive integer $j$, since the solution $u(x,t)$ is smooth, we choose $t\in[-\delta,0]$ for $0<\delta<1$, by Taylor's theorem,
\begin{equation}\label{2.14}
	u(x,t)-\sum_{i=0}^{j-1}\partial_t^iu(x,0)\frac{t^i}{i!}=\frac{t^j}{j!}\partial_s^ju(x,s),
\end{equation}
where $s=s(x,t,j)\in[t,0]$. It suffices to prove that the right hand side of \eqref{2.14} tends to zero when $j$ tends to infinity for any $x\in B_o(R)$ and $t\in [-\delta,0]$ with $\delta>0$ sufficiently small.

Since $u^2$ is a non-negative subsolution to the heat equation, we apply Proposition \ref{mean value inequality} with $p=1$. Given a point $(x_0,t_0)\in M \times [-1,0]$ and a positive integer $k$, by letting $s=t_0$, $r=\frac{1}{\sqrt{k}}$, $\tilde{\delta}=\frac{1}{2}$, $\delta'=1$ in \eqref{2.5}, we have
\begin{equation}
	u^2(x_0,t_0) 	\le\frac{c_1(n)e^{c_2(n)(K\frac{1}{k}+\sup\limits_{B_{x_0}(\frac{3}{\sqrt{k}})}|f|)}}{(\frac{1}{2})^{4n}\frac{1}{k}\vol(B_{x_0}(\frac{1}{\sqrt{k}}))}\int_{B_{x_0}(\frac{1}{\sqrt{k}})\times\left[t_0-\frac{1}{k},t_0\right]}u^2dvdt.\nonumber
\end{equation}
We observe that 
\begin{equation}
	\sup\limits_{B_{x_0}(\frac{3}{\sqrt{k}})}|f|\le\sup\limits_{B_{o}(d(x_0,o)+\frac{3}{\sqrt{k}})}|f|\le a\left(d(x_0,o)+\frac{3}{\sqrt{k}}\right)^2+b\le2ad^2(x_0,o)+\frac{18a}{k}+b,\nonumber
\end{equation}
then we have 
\begin{equation}\label{2.6}
		u^2(x_0,t_0) 	\le\frac{c_3(n)e^{c_4(n)(K+ad^2(x_0,o)+a+b)}k}{\vol(B_{x_0}(\frac{1}{\sqrt{k}}))}\int_{B_{x_0}(\frac{1}{\sqrt{k}})\times\left[t_0-\frac{1}{k},t_0\right]}u^2dvdt.
\end{equation}
Since $(\partial_t-\Delta)\partial_t^{k-1}u=0$, from \eqref{2.6}, we obtain
\begin{equation}\label{2.7}
	(\partial_t^{k-1}u)^2(x_0,t_0) \le\frac{c_3(n)e^{c_4(n)(K+ad^2(x_0,o)+a+b)}k}{\vol(B_{x_0}(\frac{1}{\sqrt{k}}))}\int_{B_{x_0}(\frac{1}{\sqrt{k}})\times\left[t_0-\frac{1}{k},t_0\right]}(\partial_t^{k-1}u)^2dvdt
\end{equation}
for a positive integer $k$.\\
Next, we will bound the right-hand side of \eqref{2.7}.\\
For positive integers $j=1,2,...,k$, we define the following domains:\\
\begin{equation}
	\Omega_j^1=B_{x_0}\left(\frac{j}{\sqrt{k}}\right)\times\left[t_0-\frac{j}{k},t_0\right],\nonumber
\end{equation}
\begin{equation}
	\Omega_j^2=B_{x_0}\left(\frac{j+0.5}{\sqrt{k}}\right)\times\left[t_0-\frac{j+0.5}{k},t_0\right].\nonumber
\end{equation}
It is easy to see that 
\begin{equation}\label{sw2}
	\Omega_j^1\subset	\Omega_j^2\subset	\Omega_{j+1}^1.
\end{equation}
Let $\psi_j^{(1)}$ be a standard Lipschitz cut-off function supported in
\begin{equation}
	B_{x_0}\left(\frac{j+0.5}{\sqrt{k}}\right)\times\left(t_0-\frac{j+0.5}{k},t_0+\frac{j+0.5}{k}\right)\nonumber
\end{equation}
  satisfying
  \begin{equation}\label{sw1}
  	\psi_j^{(1)}=1\  in\  \Omega_j^1 \  and\  |\nabla\psi_j^{(1)}|^2+|\partial_t\psi_j^{(1)}|\le Ck,
  \end{equation}
   where $C$ is a universal constant that may be changed line by line.

For the above cut-off function $\psi=\psi_j^{(1)}$, since $(\Delta-\partial_t)u=0$, using integration by parts, we compute that 
\begin{equation}
	\begin{aligned}
		\int_{\Omega_j^2}(u_t)^2\psi^2dvdt&=\int_{\Omega_j^2}u_t\Delta u \psi^2dvdt\\
		&=-\int_{\Omega_j^2} \left<(\nabla u)_t,\nabla u\right>\psi^2dvdt-\int_{\Omega_j^2}u_t\left<\nabla u,\nabla \psi^2\right>dvdt\\
		&=-\frac{1}{2}\int_{\Omega_j^2}(|\nabla u|^2)_t\psi^2dvdt-2\int_{\Omega_j^2}u_t\psi\left<\nabla u,\nabla \psi\right>dvdt\\
		&=-\frac{1}{2}\int_{B_{x_0}\left(\frac{j+0.5}{\sqrt{k}}\right)}(|\nabla u|^2\psi^2
		)(x,t_0)dv+\frac{1}{2}\int_{\Omega_j^2}|\nabla u|^2(\psi^2)_tdvdt-2\int_{\Omega_j^2}u_t\psi\left<\nabla u,\nabla \psi\right>dvdt\\
		&\le \frac{1}{2}\int_{\Omega_j^2}|\nabla u|^2(\psi^2)_tdvdt+\frac{1}{2}\int_{\Omega_j^2}(u_t)^2\psi^2dvdt+2\int_{\Omega_j^2}|\nabla u|^2|\nabla\psi|^2dvdt.\nonumber
	\end{aligned}
\end{equation}
By \eqref{sw2} and \eqref{sw1}, we have that 
\begin{equation}\label{2.8}
\int_{\Omega_j^1}(u_t)^2dvdt\le Ck\int_{\Omega_j^2}|\nabla u|^2dvdt.
\end{equation}
Let $\psi_j^{(2)}$ also be a standard Lipschitz cut-off function supported in 
\begin{equation}
	B_{x_0}\left(\frac{j+1}{\sqrt{k}}\right)\times\left(t_0-\frac{j+1}{k},t_0+\frac{j+1}{k}\right)\nonumber
\end{equation}
satisfying 
\begin{equation}\label{sw3}
	\psi_j^{(2)}=1\  in\  \Omega_{j}^2 \ and\  |\nabla\psi_j^{(2)}|^2+|\partial_t\psi_j^{(2)}|\le Ck.
\end{equation}
Then we can obtain
\begin{equation}\label{2.9}
	\int_{\Omega_j^2}|\nabla u|^2dvdt\le Ck\int_{\Omega_{j+1}^1}u^2dvdt.
\end{equation}
To achieve \eqref{2.9}, considering the cut-off function $\varphi=\psi_j^{(2)}$, by $(\Delta-\partial_t)u=0$, using integration by parts, we are continue to calculate that 
\begin{equation}
	\begin{aligned}
		&\frac{1}{2}\int_{\Omega_{j+1}^1}\partial_t(u^2\varphi^2)dvdt-\int_{\Omega_{j+1}^1}\varphi\varphi_tu^2dvdt=\int_{\Omega_{j+1}^1}uu_t\varphi^2dvdt\\
		&=\int_{\Omega_{j+1}^1}u\Delta u \varphi^2dvdt=-\int_{\Omega_{j+1}^1}|\nabla u|^2\varphi^2dvdt-2\int_{\Omega_{j+1}^1}u\varphi\left<\nabla u,\nabla \varphi\right>dvdt.\nonumber
	\end{aligned}
\end{equation}
Noticing that 
\begin{equation}
	\frac{1}{2}\int_{\Omega_{j+1}^1}\partial_t(u^2\varphi^2)dvdt=\frac{1}{2} \int_{B_{x_0}\left(\frac{j+1}{\sqrt{k}}\right)}(u^2\varphi^2)(x,t_0)dv\ge0,\nonumber
\end{equation}
we have that 
\begin{equation}
	\begin{aligned}
			\int_{\Omega_{j+1}^1}|\nabla u|^2\varphi^2dvdt&\le \int_{\Omega_{j+1}^1}\varphi\varphi_tu^2dvdt-2\int_{\Omega_{j+1}^1}u\varphi\left<\nabla u,\nabla \varphi\right>dvdt\\
			&\le \int_{\Omega_{j+1}^1}\varphi\varphi_tu^2dvdt+\frac{1}{2}\int_{\Omega_{j+1}^1}|\nabla u|^2\varphi^2dvdt+2\int_{\Omega_{j+1}^1}u^2|\nabla \varphi|^2dvdt.\nonumber
	\end{aligned}
\end{equation}
This implies that
\begin{equation}
	\int_{\Omega_{j+1}^1}|\nabla u|^2\varphi^2dvdt\le 2\int_{\Omega_{j+1}^1}\varphi\varphi_tu^2dvdt+4\int_{\Omega_{j+1}^1}u^2|\nabla \varphi|^2dvdt.\nonumber
\end{equation}
Then \eqref{2.9} follows by \eqref{sw2} and \eqref{sw3}.\\
Combining \eqref{2.8} and \eqref{2.9}, we achieve that
\begin{equation}
	\int_{\Omega_{j}^1}(u_t)^2dvdt\le Ck^2\int_{\Omega_{j+1}^1}u^2dvdt.\nonumber
\end{equation}
Since above inequality holds for all solutions of the heat equation, then we can replace $u$ by $\partial_t^ju$. By induction, we conclude that
\begin{equation}
	\int_{\Omega_{1}^1}(\partial_t^{k-1}u)^2dvdt\le C^{k-1}k^{2(k-1)}\int_{\Omega_{k}^1}u^2dvdt.\nonumber
\end{equation}
By the selection of $\Omega_{1}^1$ and $\Omega_{k}^1$, we substitute the above inequality into \eqref{2.7} to get that
\begin{equation}\label{2.10}
		(\partial_t^{k-1}u)^2(x_0,t_0) \le\frac{c_3(n)e^{c_4(n)(K+ad^2(x_0,o)+a+b)}k}{\vol(B_{x_0}(\frac{1}{\sqrt{k}}))}C^{k-1}k^{2(k-1)}\int_{\Omega_{k}^1}u^2dvdt.
\end{equation}
Using exponential quadratic  growth condition \eqref{s2} and the triangle inequality, for some point $(x,t)\in\Omega_{k}^1 $ we deduce
\begin{equation}\label{2.11}
	|u(x,t)|^2\le e^{2A_2(d(x,x_0)+d(x_0,o))^2}\le e^{2A_2(\sqrt{k}+d(x_0,o))^2}\le e^{4A_2k+4A_2d^2(x_0,o)}.
\end{equation}
By the volume comparison theorem \eqref{VC}, we have
\begin{equation}\label{2.12}
	\frac{\vol(B_{x_0}(\sqrt{k}))}{\vol(B_{x_0}(\frac{1}{\sqrt{k}}))}\le k^ne^{\frac{K}{6}\left(k-\frac{1}{k}\right)+6\sup\limits_{B_{x_0}(3\sqrt{k})}|f|}\le k^ne^{\frac{K}{6}k+6(2ad^2(x_0,o)+18ak+b)}.
\end{equation}
Substituting \eqref{2.11} and \eqref{2.12} into \eqref{2.10} gives
\begin{equation}\label{2.13}
|\partial_t^{k-1}u(x_0,t_0)|\le A_3^kk^{k-1}e^{A_4d^2(x_0,o)}
\end{equation}
for all integers $k\ge1$. Here $A_3$ and $A_4$ are two positive constants depending on $K,n,a,b,A_2$ and $n,a,A_2$, respectively.

Combining \eqref{2.13}, for \eqref{2.14}, we know that, for  $\delta<\frac{1}{A_3e}$, the right hand side of \eqref{2.14} converges to $0$ uniformly for $x\in B_o(R)$ as $j\rightarrow \infty$. Hence
\begin{equation}
	u(x,t)=\sum_{j=0}^{\infty}\partial_t^ju(x,0)\frac{t^j}{j!},\nonumber
\end{equation}
that is, $u(x,t)$ is time analytic with radius $\delta$. Write $a_j=a_j(x)=\partial_t^ju(x,0)$. We have that 
\begin{equation}
	\partial_t u(x,t)=\sum_{j=0}^{\infty}a_{j+1}(x)\frac{t^j}{j!} \ and\ \Delta u(x,t)=\sum_{j=0}^{\infty}\Delta a_j(x)\frac{t^j}{j!},\nonumber
\end{equation}
where both series converge uniformly for $(x,t)\in B_o(R)\times[-\delta,0]$. Since $(\Delta-\partial_t)u=0$, this gives that 
\begin{equation}
	\Delta a_j(x)=a_{j+1}(x)\nonumber
\end{equation}
and 
\begin{equation}
	|a_j(x)|\le A_3^{j+1}(j+1)^je^{A_4d^2(x,o)}.\nonumber
\end{equation}
Here $A_3$ and $A_4$ are two positive constants depending on $K,n,a,b,A_2$ and $n,a,A_2$, respectively.\\ \qed

Next we apply Theorem \ref{1.1} to prove Corollary \ref{1.2}.

\noindent{\it Proof of Corollary \ref{1.2}.} Assume that $u(x,t)$ is a smooth solution to  \eqref{4} with exponential quadratic  growth of the space variable in $M\times(0,\delta)$. Then 
\begin{equation}
	(\Delta-\partial_t)u(x,-t)=0\ and\ |u(x,-t)|\le A_1e^{A_2d^2(x,o)},\nonumber
\end{equation}
where $A_1$ and $A_2$ are some positive constants.
By Theorem \ref{1.1}, we have 
\begin{equation}
	u(x,-t)=\sum_{j=0}^{\infty}a_j(x)\frac{(-t)^j}{j!}.\nonumber
\end{equation}
Combining the initial condition of \eqref{4} with  Theorem \ref{1.1}, then \eqref{5} follows.

 On the other hand, suppose \eqref{5} holds. Setting $u(x,t)=\sum_{j=0}^{\infty}\Delta^ja(x)\frac{t^j}{j!}$, by \eqref{5}, it is easy to see that 
\begin{equation}
	\sum_{j=0}^{\infty}\Delta^{j+1}a(x)\frac{t^j}{j!}\ and\ \sum_{j=0}^{\infty}\Delta^{j}a(x)\frac{\partial_tt^j}{j!}\nonumber
\end{equation}
all converge absolutely and uniformly in $B_o(R)\times[-\delta,0]$ for any fixed $R>0$ and $\delta>0$ sufficiently small. Hence
\begin{equation}
	\left(\Delta-\partial_t\right)u(x,t)=0\ for\ (x,t)\in M\times[-\delta,0]. \nonumber
\end{equation}
By \eqref{5} again, we get the exponential quadratic growth for $u$,
\begin{equation}
	|u(x,t)|\le \sum_{j=0}^{\infty}|\Delta^ja(x)|\frac{|t|^j}{j!}\le A_3e^{A_4d^2(x,o)}\sum_{j=0}^{\infty}\frac{(A_3(j+1)|t|)^j}{j!}\le A_5e^{A_4d^2(x,o)}\nonumber
\end{equation}
provided that $t\in[-\delta,0]$ with  $\delta>0$ sufficiently small. 

Then a smooth solution with desired growth condition to \eqref{4} follows by letting $u=u(x,-t)$. \\ \qed

The proof of Theorem \ref{1.3} is similar to Theorem \ref{1.1}. We only present the key steps.
\noindent{\it Proof of Theorem \ref{1.3}.}  From \eqref{2.10}, we know, for $(x_0,t_0)\in M\times[-1,0]$ and any positive integer $k$,
\begin{equation}\label{2.15}
\begin{aligned}
		(\partial_t^{k-1}u)^2(x_0,t_0) \le\frac{\bar{c}_1(n)e^{\bar{c}_2(n)(K+ad^2(x_0,o)+a+b)}k}{\vol(B_{x_0}(\frac{1}{\sqrt{k}}))}C^{k-1}k^{2(k-1)}\int_{B_{x_0}(\sqrt{k})\times[t_0-1,t_0]}u^2dvdt.
\end{aligned}
\end{equation} 
By mean value theorem, there exists $\xi\in (t_0-1,t_0)$ such that 
\begin{equation}\label{2.18}
	\begin{aligned}
		\int_{B_{x_0}(\sqrt{k})\times[t_0-1,t_0]}u^2dvdt&=\int_{B_{x_0}(\sqrt{k})}u^2(x,\xi)dv\\
		&\le\left(\int_{B_{x_0}(\sqrt{k})}|u|^p(x,\xi)dv\right)^{\frac{2}{p}}\vol(B_{x_0}(\sqrt{k}))^{1-\frac{2}{p}}\\
		&\le L^2\vol(B_{x_0}(\sqrt{k}))^{1-\frac{2}{p}}
	\end{aligned}
\end{equation}
for $p\ge2$, where we used H\"older inequality in the second line and in the last line we used the assumption \eqref{ssss4}.\\
By volume comparison theorem \eqref{VC} and $k\ge1$, we have 
\begin{equation}\label{2.19}
	\begin{aligned}
	\frac{\vol(B_{x_0}(\sqrt{k}))^{1-\frac{2}{p}}}{\vol(B_{x_0}(\frac{1}{\sqrt{k}}))}&=\frac{\vol(B_{x_0}(\sqrt{k}))\vol(B_{x_0}(\sqrt{k}))^{-\frac{2}{p}}}{\vol(B_{x_0}(\frac{1}{\sqrt{k}}))}\\
	&\le(k)^ne^{\frac{K}{6}k+12ad^2(x_0,o)+108ak+6b}\vol(B_{x_0}(1))^{-\frac{2}{p}}.
	\end{aligned}
\end{equation}
To get a lower bound of $\vol(B_{x_0}(1))$, we use the volume comparison theorem \eqref{VC} again, then 
\begin{equation}\label{2.20}
	\begin{aligned}
			\vol(B_{o}(1))&\le \vol (B_{x_0}(d(x_0,o)+1))\\
			&\le \vol(B_{x_0}(1))(d(x_0,o)+1)^ne^{\frac{K}{6}[(d(x_0,o)+1)^2-1]+6\sup\limits_{B_{x_0}(3d(x_0,o)+3)}|f|}\\
			&\le\vol(B_{x_0}(1))e^{nd(x_0,o)+\frac{K}{6}(d^2(x_0,o)+2d(x_0,o))+6\sup\limits_{B_{o}(4d(x_0,o)+3)}|f|}\\
			&\le\vol(B_{x_0}(1))e^{(n+\frac{K}{2})d^2(x_0,o)+\frac{1}{4}n+\frac{K}{12}+6(32ad^2(x_0,o)+18a+b)}.
	\end{aligned}
\end{equation}
Combining \eqref{2.20}, \eqref{2.19}, \eqref{2.18} with \eqref{2.15}, we arrive at 
\begin{equation}
	|\partial_t^{k-1}u(x_0,t_0)|\le A_6^kk^{k-1}e^{A_7d^2(x_0,o)}\vol(B_o(1))^{-\frac{1}{p}}L\nonumber
\end{equation}
for all integers $k\ge1$. Here $A_6$ and $A_7$ are two positive constants depending on $n,K,a,b,p$ and $n,a,K,p$, respectively. \\ \qed
\section*{Acknowledgements}

Research is partially supported by NSFC Grant No. 11971168, Shanghai Science and Technology Innovation Program Basic Research Project STCSM 20JC1412900, and Science and Technology Commission of Shanghai Municipality (STCSM) No. 22DZ2229014.
\section*{Data availibility}
The datasets generated during and/or analysed during the current study are available from the corresponding author on reasonable request.


\begin{thebibliography}{99}
	\bibitem{BE}  D. Bakry; M. \'Emery, \emph{Diffusions hypercontractives.}  (French) [Hypercontractive diffusions] S\'eminaire de probabilit\'es, XIX, 1983/84, 177-206, Lecture Notes in Math., 1123, Springer, Berlin, 1985.	
	\bibitem{Chen} Bing-Long Chen, Strong uniqueness of the Ricci flow, J. Differential Geom. 82(2): 363-382
	\bibitem{Gi} Giga Y. \emph{Time and spatial analyticity of solutions of the Navier-Stokes equations}. Comm Partial Differential Equations, 1983, 8(8): 929–948
	\bibitem{Ham} R. Hamilton, {\it The formation of singularities in the Ricci flow}, Surveys in differential geometry, Vol. II (Cambridge, MA, 1993), 7136, Int. Press, Cambridge, MA, 1995.
	\bibitem{LV}  John Lott; C\'edric Villani, \emph{Ricci curvature for metric-measure spaces via optimal transport.} Ann. of Math. (2) 169 (2009), no. 3, 903-991.
	\bibitem{LSC} Escauriaza L, Montaner S, Zhang C. \emph{Analyticity of solutions to parabolic evolutions and applications}. SIAM J Math Anal, 2017, 49(5): 4064–4092
	\bibitem{Per} G. Perelman, \emph{The entropy formula for the Ricci flow and its geometric applications,} arXiv:math/0211159.
	\bibitem{PRS} S. Pigola, M. Rimoldi, A. Setti, {\it Remarks on non-compact gradient Ricci solitons}, Math. Z. 268 (2011), 777-790.
	\bibitem{St1}  Karl-Theodor Sturm, \emph{On the geometry of metric measure spaces.I.} Acta Math. 196 (2006), no. 1, 65-131.
	\bibitem{St2}  Karl-Theodor Sturm, \emph{On the geometry of metric measure spaces. II.} Acta Math. 196 (2006), no. 1, 133-177.
	\bibitem{SWZ} XingYu Song; Ling Wu; Meng Zhu, \emph{Heat kernel estimate for the Laplace-Beltrami operator under Bakry-\'Emery Ricci curvature condition and applications}. submitted.
    \bibitem{WJY}Wu, Jiayong, \emph{Time analyticity for the heat equation on gradient shrinking Ricci solitons}. Acta Math. Sci. Ser. B (Engl. Ed.) 42 (2022), no. 4, 1690–1700.
    \bibitem{WLL} Han F W, Hua B B, Wang L L. \emph{Time analyticity of solutions to the heat equation on graphs}. Proc Amer Math Soc, 2021, 149(6): 2279–2290
    \bibitem{ZQ2}Dong, Hongjie ; Zhang, Qi S., \emph{Time analyticity for the heat equation and Navier-Stokes equations}. J. Funct. Anal. 279 (2020), no. 4, 108563, 15 pp.
	\bibitem{ZQ1}Zhang, Qi S., \emph{A note on time analyticity for ancient solutions of the heat equation}. Proc. Amer. Math. Soc. 148 (2020), no. 4, 1665–1670.
	\bibitem{ZS}S. Zhang, \emph{On a sharp volume estimate for gradient Ricci solitons with scalar curvature bounded below}, Acta Math. Sin. 27(5) (2011), 871–882.
    \bibitem{WD}Widder, D. V, \emph{Analytic solutions of the heat equation}. Duke Math. J. 29 (1962), 497–503.
\end{thebibliography}
\end{document}